\documentclass[12pt]{article}
\usepackage[english]{babel}
\usepackage{amssymb,amsmath,graphics}
\parskip\medskipamount
\newtheorem{theorem}{Theorem}[section]
\newtheorem{Theorem}{Main Result}
\newtheorem{remark}[theorem]{Remark}

\newtheorem{lemma}[theorem]{Lemma}
\newtheorem{cor}[theorem]{Corollary}
\def\<{\langle}
\def\>{\rangle}
\newcommand{\proof}{\emph{Proof.~}}

\newcommand{\PG}{\mathsf{PG}}

\newcommand{\cP}{\mathcal{P}}

\newcommand{\cM}{\mathcal{M}}

\newcommand{\sfF}{\mathsf{F}}

\def\qed{{\hfill\hphantom{.}\nobreak\hfill$\Box$}}

\newcommand{\cH}{\mathcal{H}}

\newcommand{\ssW}{\mathsf{W}}
\begin{document}

\author{Koen Struyve\thanks{The author is supported by  the Fund for Scientific Research --
Flanders (FWO - Vlaanderen)} }
\title{\bf Quadrangles embedded in metasymplectic spaces}

\maketitle

\begin{abstract}
During the final steps in the classification of the Moufang quadrangles by Jacques Tits and Richard Weiss a new class of Moufang quadrangles unexpectedly turned up. Subsequently Bernhard M\"uhlherr and Hendrik Van Maldeghem showed that this class arises as the fixed points and hyperlines of certain involutions of a metasymplectic space (or equivalently a building of type $\sfF_4$). In the same paper they also showed that other types of Moufang quadrangles can be embedded in a metasymplectic space as points and hyperlines.

In this paper, we reverse the question: given a (thick) quadrangle embedded in a metasymplectic space as points and hyperlines, when is such a quadrangle a Moufang quadrangle?
\end{abstract}

\section{Introduction} 
Generalized polygons are the geometries related to the spherical rank 2 buildings. These geometries were introduced by Jacques Tits in the appendix of \cite{Tit:59} prior to the first formal definition of buildings in the literature. The first examples of generalized polygons mainly arose as \emph{embeddings} in projective spaces, i.e., the points of the polygon were some points of a projective space, while the lines of the polygon could be identified with some lines of the projective space, with natural incidence relation.   If the embedding is `nice', then it automatically inherits beautiful symmetry properties from the projective space; see \cite{Die:80,Lef:81,Ste-Mal:00,Ste-Mal:04,Tha-Mal:04}. `Nice' could mean that the lines of the polygon through any point are contained in a certain subspace of the projective space (plane, hyperplane), or that the the points not opposite a given point in the polygon do not span the entire projective space, or just a bound on the dimension of the projective space together with the fact that all points of the projective space on any line of the polygon belong to the polygon. In particular, the previous references contain characterizations and classifications of the `nice' embeddings of the Moufang generalized quadrangles and hexagons. (The `Moufang' condition is a condition on the automorphism group of the polygon, implying a lot of symmetries. The Moufang polygons are classified in \cite{titsweiss}.) 

However, not all Moufang polygons admit an embedding as considered above. The notable examples are the exceptional Moufang quadrangles and their duals, the duals of some embeddable classical Moufang quadrangles, and the duals of the exceptional Moufang hexagons and of the Ree-Tits octagons. These exceptional polygons geometrically arise in a different way: they do not arise from `forms' of a projective space, but from `forms' of buildings of exceptional type and rank at least 4. All types arise: $\mathsf{E}_6,\mathsf{E}_7,\mathsf{E}_8,\sfF_4$.   In this paper, we take a closer look at the situation of $\sfF_4$ (or called metasymplectic spaces from a geometric point of view). This case is the least `algebraic' of the lot. Indeed, buildings of type $\sfF_4$ give rise to octagons and to quadrangles, but the corresponding `forms' are not forms of an algebraic group. Instead, they owe their existence to the exceptional behaviour of fields of characteristic 2, and the related existence of groups and buildings of `mixed' type, see \cite{Tit:74}. The situation of the mixed case being somewhat less algebraic means also that it is somewhat more geometric.   This is the starting point of the present paper. Our goal is to find a `nice' property of the embedding of the exceptional Moufang quadrangles in buildings of type $\sfF_4$ that guarantees that \emph{any} quadrangle embedded in a building of type $\sfF_4$ with that property, is automatically a Moufang quadrangle. This property will be denoted by (OV) below. Roughly, we require that the points of the quadrangle are points of the building, the lines of the quadrangle are hyperlines of the building (with natural incidence), and (OV) says that any two noncollinear points of the quadrangle are never contained in a hyperline of the building. In other words, \emph{collinearity in the quadrangle coincides with cohyperlinearity in the building}. This very natural property surprisingly is enough to characterize the Moufang quadrangles arising from buildings of type $\sfF_4$.      

\section{Preliminaries and Main Results}
\subsection{Generalized polygons and metasymplectic spaces}
A \emph{flag} of a geometry is a set of mutually incident elements, a \emph{chamber} is a maximal flag.

A \emph{generalized $n$-gon} is a geometry with two types of elements, \emph{points} and \emph{lines} and a (symmetric) incidence relation, such that each two elements are contained in an ordinary $n$-gon, there are no ordinary $k$-gons with $2 \leq k < n$, and such that there exists an ordinary $n+1$-gon. This last condition assures that the generalized $n$-gon is \emph{thick}: there are at least three points on a line and three lines through a point.

The flag complexes of these geometries form the spherical buildings of rank 2, the case of a generalized triangle is better known as a (axiomatic) projective plane. An  \emph{apartment} of a generalized $n$-gon is an ordinary $n$-gon. Often we will omit `generalized' if the context is clear. 

A special class of polygons are the \emph{Moufang polygons}. These are polygons with extra group-theoretic conditions. We only note that these have been classified by Jacques Tits and Richard Weiss in 2002 (\cite{titsweiss}).

We use the following definition for metasymplectic spaces (\cite[p. 79]{hvm1}): a \emph{metasymplectic space} $\cM$ is a geometry with four types of elements, denoted with \emph{points}, \emph{lines}, \emph{planes} and \emph{hyperlines} and a (symmetric) incidence relation satisfying the four axioms listed below.

A \emph{residue} of a flag $A$ is the geometry of elements distinct to those of $A$ and incident with all elements of $A$. The \emph{type} of a flag $A$ is the set of types of its elements.
\begin{itemize}
\item[{(M1)}]
The residue of any flag of type \{point, line\} or \{plane, hyperline\} is a projective plane.
\item[{(M2)}]
The residue of any flag of type \{point, plane\}, \{line, hyperline\} or \{line, plane\} is a generalized digon.
\item[{(M3)}]
The residue of any flag of type \{point, hyperline\} is a generalized quadrangle.
\item[{(M4)}]
Two distinct non-point elements have different sets of points incident with them.
\end{itemize}
Using (M1) to (M4), one can show the dual property of (M4), making the definition self-dual. The flag complexes of these metasymplectic spaces form the buildings of type $\sfF_4$.

\begin{remark}
\rm Instead of the notion hyperline, some authors use the term `symplecton'.
\end{remark}

\subsection{Embeddings of quadrangles in the metasymplectic space}
We consider embeddings of the following kind: given a metasymplectic space $\cM$ together with a set $\cP$ of points of $\cM$, and a set $\cH$ of hyperlines of $\cM$ such that the incidence relation defined on them by taking the restriction of the incidence relation of $\cM$ defines a generalized quadrangle $\Gamma$. We then say that the quadrangle $\Gamma$ is \emph{point--hyperline embedded} in $\cM$.

Examples of such embeddings are constructed by Hendrik Van Maldeghem en Bernhard M\"uhlherr in~\cite{hvmmhl}. There it is shown that the exceptional Moufang quadrangles of type $\sfF_4$ and certain mixed quadrangles appear as fixed point structures of involutions of metasymplectic spaces over fields with characteristic 2. As the subquadrangles of a point--line embedded quadrangle will also be point--line embedded, orthogonal and symplectic quadrangles also appear. All these quadrangles are Moufang and share the property that no two points of the quadrangle are collinear in the metasymplectic space.

Embeddings will be denoted \emph{improper} if each two hyperlines in $\cH$ incident with a point in $\cP$ always share a line. If this is the case, then it can be shown that all the hyperlines in $\cH$ incident with a certain point $p$ in $\cP$ share a line $L_p$. By substituting each point with its associated line, it follows that we can view the quadrangle embedded as lines and hyperlines. 

We now construct an example of an improper embedding. Let $\{p,L\}$ be an incident point--line pair of a metasymplectic space $\cM$ which is defined over some field containing the finite field over four elements.  The residue forms a projective plane, containing a sub-projective plane isomorphic to $\PG (2,4)$. The symplectic quadrangle $\ssW(2)$ can be embedded in this plane (see~\cite{Beu:86}). Returning to our metasymplectic space $\cM$, we thus have $\ssW(2)$ embedded in $\cM$ as planes and hyperlines. Now choose for each plane of this embedding a point incident with the plane, which produces a point--hyperline embedding. If the field which defines the metasymplectic space is `large enough', it is clear that the choices can be made such that no two collinear points of the quadrangle are collinear in the metasymplectic space. 

\begin{remark}
\rm All of the known embeddings such that no two points of the quadrangle are collinear in the metasymplectic space occur in characteristic 2 or are improper. The existence of the known proper embeddings originates from an algebraic setting, however this algebraic setting does not yield such embeddings for odd characteristic. It thus can be conjectured that these only occur in characteristic 2. More about the underlying algebraic setting can be found in~\cite[App. C]{hvm1}.
\end{remark}

\subsection{Main Result}
We now pose the inverse question: when is a point--hyperline embedded quadrangle Moufang?

\begin{Theorem}Let $\Gamma$ be a generalized quadrangle point--hyperline embedded in a metasymplectic space $\cM$, with $\cP$ the set of points and $\cH$ the set of lines of the quadrangle. Then $\Gamma$ will be a Moufang quadrangle or improperly embedded if the following property holds:
\begin{itemize}
\item[{(OV)}] No 2 points of $\cP$ in the same hyperline of $\cH$ are collinear in $\cM$.
\end{itemize}
\end{Theorem}

\begin{remark}
\rm Note that our definition of generalized polygon asks that $\Gamma$ is thick, if this would not be the case, then counter examples occur.
\end{remark}

\begin{remark}
\rm It can be shown that the residue of a hyperline forms a polar space (see property (M9) in the next section). Condition (OV) then reformulates to: the points of $\cP$ in the same hyperline of $\cH$ form a partial ovoid of the corresponding polar space.
\end{remark}

\section{Proof of the main result}
Suppose we have a $\cM,\Gamma,\cP,\cH$ as given in the statement of the main result, but we do not require that the property (OV) holds yet.

If we refer to a point or line, we mean a point or line of the metasymplectic space, unless explicitly noted otherwise.

\subsection{Further concepts and some lemmas about metasymplectic spaces}
The following lemma can be found in~\cite[p. 80]{hvm1}, we will not reproduce the proof here.
\begin{lemma} 
\begin{itemize}
\item[{(M5)}] Let $x$ and $y$ be two points of $\cM$. Then one of the following situations occurs:
\begin{itemize}
\item
$x=y$.
\item
There is a unique line incident with both $x$ and $y$. In this case, we call $x$ and $y$ \emph{collinear}.
\item
There is a unique hyperline incident with both $x$ and $y$. In this case there is no line incident with both $x$ and $y$, and we call $x$ and $y$ \emph{cohyperlinear}.
\item
There is a unique point $z$ collinear with both $x$ and $y$. In this case we call $x$ and $y$ \emph{almost opposite}. In this case we call $x$ and $y$ opposite.
\item
There is no point collinear with both $x$ and $y$.
\end{itemize}
\item[{(M6)}] The intersection of two hyperlines is either empty, or a point, or a plane.
\item[{(M7)}] Let $x$ be a point and $h$ a hyperline of $\cM$. Then one of the following situations occurs:
\begin{itemize}
\item $x\in h$.
\item
There is a unique line $L$ in $h$ such that $x$ is collinear with all points of $L$. Every point $y$ of $h$ which is collinear with all points of $L$ is cohyperlinear with $x$ and the unique hyperline containing both also contains $L$. Every other point $z$ of $h$ is almost opposite $x$ and the unique point collinear with both lies on $L$.
\item
There is a unique point $u$ of $h$ cohyperlinear with $x$, the hyperline containing $x$ and $u$ only has $u$ in common with $h$. All points $v$ of $h$ collinear with $u$ are almost opposite $x$, and the point collinear with both does not lie in $h$. All points $w$ of $h$ cohyperlinear with $u$ are opposite $x$.
\end{itemize}
\item[{(M9)}] The residue of a hyperline forms a polar space.
\qed
\end{itemize} 

\end{lemma}
Again note that the dual statements also hold. Property (M8) given in~\cite{hvm1} is omitted as we will not need it here. 

A \emph{path of chambers} is an ordered set of chambers such that each chamber only differs one element from the previous. The number of chambers in the path minus one is the \emph{length} of the path. Let $W$ be the spherical Coxeter group of type $F_4$, this is the group generated by symbols $s_1,s_2,s_3, s_4$ and identity element $e$, with relations $(s_is_j)^{m_{ij}} = e$, with $m_{ij}$ as given by the following matrix:
\begin{equation*}
(m_{ij})=\left(\begin{array}{cccc}1 & 3 & 2 & 2 \\3 & 1 & 4 & 2 \\2 & 4 & 1 & 3 \\2 &2 & 3 & 1\end{array}\right)
\end{equation*} 
With a path one can associate a word in the symbols $s_1,s_2,s_3$ and $s_4$. For each chamber beyond the first one sets $s_1,s_2,s_3$ or $s_4$ respectively if that chamber and the previous chamber differ in a point, line, plane or hyperline respectively. The following lemma is well known in the theory of buildings (see for example in~\cite{ronan}):
\begin{lemma}
A path between two chambers has the shortest length possible between those two chambers, if and only if the associated word has no shorter representation in the Coxeter group $W$.
\qed
\end{lemma} 
The spherical Coxeter group $W_{\{1,2,3\}}$ is the subgroup of $W$ generated by $s_1, s_2$ and $s_3$, analogously $W_{\{2,3,4\}}$ will be the subgroup generated by $s_2$, $s_3$ and $s_4$. 

\begin{lemma} \label{lemma:red}
The following double cosets are written in such a way that the representative is of shortest length:
\begin{itemize}
\item $W_{\{2,3,4\}} s_1s_2s_3s_2s_1 W_{\{2,3,4\}}, W_{\{1,2,3\}} s_4s_3s_2s_3s_4 W_{\{1,2,3\}}$
\item $ W_{\{1,2,3\}} s_4s_3s_2s_3s_4s_1s_2s_3s_2s_1 W_{\{2,3,4\}}, W_{\{2,3,4\}} s_1s_2s_3s_2s_1s_4s_3s_2s_3s_4  W_{\{1,2,3\}} $
\item $ W_{\{2,3,4\}} s_1s_2s_3s_2s_1s_4s_3s_2s_3s_4s_1s_2s_3s_2s_1 W_{\{2,3,4\}}$, \\ $W_{\{1,2,3\}} s_4s_3s_2s_3s_4s_1s_2s_3s_2s_1s_4s_3s_2s_3s_4 W_{\{1,2,3\}}$
\item $ W_{\{2,3,4\}} s_1s_2s_3s_2s_1s_4s_3s_2s_3s_4s_1s_2s_3s_2s_1s_4s_3s_2s_3s_4 W_{\{1,2,3\}}$, \\ $W_{\{1,2,3\}} s_4s_3s_2s_3s_4s_1s_2s_3s_2s_1s_4s_3s_2s_3s_4s_1s_2s_3s_2s_1 W_{\{2,3,4\}}$
\end{itemize}
\end{lemma}
\proof
Lengthy but straightforward calculations. \qed

Given two flags $A$ and $B$ of the metasymplectic space $\cM$. Consider the shortest paths from chambers containing $A$ to chambers containing $B$ (minimized over all choices of such chambers). The intersection of the last chambers in all these paths is called the \emph{projection} of the flag $A$ onto $B$. A set of elements of elements of $\cM$ is called \emph{convex} if the projection between two flags of elements in that set, is a subset of that set. The following important theorem by Bernhard M\"uhlherr and Hendrik Van Maldeghem (\cite{Mhl-Mal:02}) gives us more information about convex subbuildings (as with generalized polygons, all buildings considered are thick).
\begin{theorem}
A convex subbuilding of a Moufang building is again a Moufang building.
\qed
\end{theorem}
Or applied to our case (buildings of type $\sfF_4$ are always Moufang):
\begin{cor}\label{cor:convex}
A convex point--hyperline embedded quadrangle $\Gamma$ in a metasymplectic space $\cM$ is Moufang.
\qed
\end{cor}
\subsection{Embedding apartments}
First we investigate how the apartments of the quadrangle are embedded in $\cM$. Let $\{p,h\}, \{q,g\}$ $(p,q \in \cP ,h,g \in \cH)$ be 2 chambers of $\Gamma$ such that $p \notin g,q \notin h$ and the hyperlines $h$ and $g$ intersect in a point or plane (these are the only possibilities barring equal or disjoint hyperlines due to (M6)). Collinearity and opposition will be used in relative to the metasymplectic space $\cM$ and not the quadrangle $\Gamma$, unless stated otherwise. 
\begin{lemma}
If $h$ and $g$ intersect in a point $u$, then one of the following holds:
\begin{itemize}
\item The points $p$ and $q$ are opposite and both are hypercollinear with $u$.
\item The points $p$ and $q$ are almost opposite and at least one point is collinear with $u$.
\item The points $p$ and $q$ are cohyperlinear and both are collinear with $u$.
\item The points $p$ and $q$ are collinear and both are collinear with $u$.
\end{itemize}
\end{lemma}
\proof
\begin{itemize}
\item
If $p$ and $q$ are opposite then (M7) applied to the point $p$ and hyperline $g$ tells us that there is exactly one point of $g$ cohyperlinear with $p$, therefore $u$ will be this point. It now follows that $p$ and $q$ both are cohyperlinear with $u$.
\item 
If $p$ and $q$ are almost opposite then applying (M7) to $p$ and $g$ leaves us with two possibilities. If there is a unique point (this point will again be $u$) of $g$ cohyperlinear with $p$ then $q$ will be collinear with $u$. If on the other hand there is a unique line $L$ in $g$ of points collinear with $p$ then the possibility that $u$ is cohyperlinear with $p$ implies that $u$ is collinear with all points of $L$ and that $h$ contains $L$. But $h$ and $g$ intersect in a point and not in line, so $p$ is collinear with $u$.
\item 
If $p$ and $q$ are cohyperlinear then again applying (M7) to $p$ and $g$ gives us that there is a line $L$ in $g$ of points collinear with $p$ (the other possibility for cohyperlinearity would imply that $u = q$ which is ruled out). If $u$ would be cohyperlinear with $p$ then $h$ and $g$ would intersect in a line as explained in the previous point, so $p$ is collinear with $u$. Interchanging the roles of $p$ and $q$ gives that both points are collinear with $u$.
\item
In the last case where $p$ is collinear with $q$, property (M7) implies that $p$ is collinear with all the points of a line $L$ of $g$. If $u$ would be cohyperlinear with $p$ then the unique hyperline $h$ containing $u$ and $p$ would also contain $q$, which is impossible. It follows that $u$ is collinear $p$ and also with $q$ which is proven analogously.\qed
\end{itemize}

\begin{lemma}
If $h$ and $g$ intersect in a plane $\pi$, then $p$ and $q$ are not opposite.
\end{lemma}
\proof
If this was the case then $p$ and $q$ would be on distance 3 from each other, but (M9) gives us that the points on distance 1 from $p$ in $\pi$ will be the points on a line of $\pi$, the same holds for $q$. Two 2 lines in a plane always have at least one point in common, thus the distance between $p$ and $q$ would be 2, producing a contradiction. \qed

Given an apartment in $\Gamma$ consisting of the points $p,q,r,s$ and hyperlines denoted by $pq,qr,rs,sp$. If the points $p$ and $r$ are opposite then the two lemmas above imply that if two points of the apartment are collinear in $\Gamma$, they are cohyperlinear in $\cM$. The hyperline $pq$ intersects $qr$ in a point, the same holds for $sp$ and $rs$. The other mutual positions divide in 2 possibilities due to the third lemma :

\begin{itemize}
\item
The hyperlines $pq$ and $sp$ intersect in a point. Then $q$ and $s$ are opposite and $qr$ and $rs$ also intersect in a point.
\item 
The hyperlines $pq$ and $sp$ intersect in a plane. Then $q$ and $s$ are not opposite and $qr$ and $rs$ also intersect in a plane.
\end{itemize}
We now state a lemma which will be used to 'reduce' the quadrangle. 
\begin{lemma}\label{LM:R}
If a set $X$ consists of mutually collinear points of $\cM$ , then this set is contained in a plane.
\end{lemma}
\proof
Let $x\in X$ be a point, taking the residue of this point, we obtain a dual rank 3 polar space where the lines $xy$ with $y\in X \backslash \{x\}$ form dual generators. All these generators intersect in lines of the polar space. If we have a proper 'triangle' of these generators and lines, and the lines meet in a single point, taking the residue again of this point, we would have a proper triangle in a quadrangle, which is impossible. So all the generators $xy$ with $y\in X \backslash \{x\}$ share at least one line, translating this back to $\cM$ we obtain that all points are contained in a plane. \qed

\subsection{Embedding quadrangles}
\subsubsection{Condition (OV)}
From now on suppose that condition (OV) holds. Let $\Sigma$ be an apartment of $\Gamma$. If 2 hyperlines of $\Sigma$ which intersect in $\Gamma$ share a point, then there occurs an opposite pair of points (opposite in $\cM$), thus the other 2 hyperlines in $\Sigma$ must also intersect in a point according to the previous section. Because the projectivity group of a point of our quadrangle is 2-transitive on the (hyper)lines through that point, either any 2 hyperlines in $\cH$ which intersect in $\Gamma$ share a point, or all hyperlines in $\cH$ which intersect in $\Gamma$ share a plane.

In the second case we can replace each point $p \in \cP$ with a line $L_p$ such that all hyperlines of $\cH$ through $p$ contain that line (this is possible due to the dual of lemma~\ref{LM:R}), so we obtain a quadrangle consisting of lines and hyperlines where no 2 lines which are collinear in the quadrangle are contained in one plane (otherwise the points corresponding to the 2 lines would be collinear in $\cM$), so we are in the improper case.

In the first case we have that 2 points of $\cP$ are cohyperlinear if they are collinear in $\Gamma$ and opposite if they are not. For hyperlines in $\cH$ we have the dual properties. In the next section we will show convexity of quadrangles within $\cM$ with such properties.

\subsubsection{Convexity of quadrangles}
In this section we prove that if we have that 2 points of $\cP$ are cohyperlinear if they are collinear in $\Gamma$ and opposite if they are not and the dual properties for hyperlines in $\cH$, then the embedded quadrangle $\Gamma$ is convex in $\cM$.


The next lemma gives us the needed building blocks for the rest of the proof of convexity.
\begin{lemma}
Let $h$ be a hyperline and $p, q$ be 2 cohyperlinear points in $h$. If we have a chamber $C$ containing $p$ and $h$, then there is a shortest path with associated word $s_1s_2s_3s_2s_1$ from $C$ to a chamber containing $q$.
\end{lemma}
\proof
The residue of $h$ will be a rank 3 polar space with $p$ and $q$ opposite points in it. The theory of buildings tells us that we can embed the flags $C \backslash \{h \}$ and $\{q\}$ of this polar space in an octahedron (this forms an apartment of the rank 3 polar space, see~\cite{ronan}). In this octahedron it is easily seen that there is a shortest path with associated word $s_1s_2s_3s_2s_1$ from $C$ to a chamber containing $q$ and $h$. Because this word is a shortest representation of the corresponding element in the group $W$, this will be a shortest path. \qed



Let $A$ and $B$ be two flags of $\Gamma$. It is clear that there exists a shortest path $\gamma_\Gamma$ in $\Gamma$ between these 2 flags starting from a chamber $C$ in $\Gamma$ containing $A$ to a chamber $D$ containing $B$. Using the above lemma (and the dual statement) to 'lift' this path to a path $\gamma_\cM$ in $\cM$, we obtain paths from each chamber containing $C$ (now viewed as flags in $\cM$) to a certain chamber containing $D$ (viewed as flag in $\cM$) with words consisting of an alternating consecution of the `building block' $s_1s_2s_3s_2s_1$ and the dual $s_4s_3s_2s_3s_4$. The Lemma~\ref{lemma:red} implies that these are also shortest paths between chambers containing $A$ and chambers containing $B$ in $\cM$. Because the paths can start from each chamber containing $C$, the projection of $B$ to $A$ will thus be completely contained within $C$ and thus within the subbuilding $\Gamma$; hence the embedded quadrangle $\Gamma$ is convex. Corollary~\ref{cor:convex} now implies that the quadrangle $\Gamma$ is Moufang.


\end{document}